%% file: ks.tex
\documentclass[11pt]{article}
\usepackage{amsmath,amsfonts,amsthm,amssymb,color,newclude}

\input{defs}

\begin{document}

\title{
Interlacing Families II: 
Mixed Characteristic Polynomials and The Kadison-Singer Problem
\thanks{
This research was partially supported by NSF grants CCF-0915487 and  CCF-1111257,
  an NSF Mathematical Sciences Postdoctoral Research Fellowship, Grant No. DMS-0902962,
  a Simons Investigator Award to Daniel Spielman, and a MacArthur Fellowship.
}}

\author{
Adam W. Marcus\\
Yale University\\
\and
Daniel A. Spielman \\ 
Yale University
\and 
Nikhil Srivastava\\
Microsoft Research, 
India
}

\maketitle

\begin{abstract}
We use the method of interlacing families of polynomials introduced in \cite{lifts_paper} to prove two theorems known to
  imply a positive solution to the Kadison--Singer problem.
  
The first is Weaver's conjecture $KS_{2}$ \cite{weaver}, which is known to   
  imply Kadison--Singer via a projection paving conjecture
  of Akemann and Anderson \cite{akemann1991lyapunov}.
  The second is a formulation due to Casazza, et al., \cite{CEKP07} of Anderson's original paving conjecture(s) \cite{anderson1979a,anderson1979b,anderson1981}, for which we are able to compute explicit paving bounds.
    
The proof involves an analysis of the largest roots of a family of polynomials that we call the
  ``mixed characteristic polynomials''
  of a collection of matrices.
\end{abstract}

\newpage

\include*{ksintro}

\include*{kspreliminaries}

\include*{ksmixed}
\include*{ksbarriers}

\include*{kspaving}

\include*{ksconclusion}

\bibliographystyle{dansalpha}
\bibliography{ks}

\end{document}

%% file: defs.tex
\newtheorem{theorem}{Theorem}[section]
\newtheorem{corollary}[theorem]{Corollary}
\newtheorem{lemma}[theorem]{Lemma}

\newtheorem{proposition}[theorem]{Proposition}

\newtheorem{conjecture}[theorem]{Conjecture}

\newtheorem{question}[theorem]{Question}

\theoremstyle{definition}
\newtheorem{definition}[theorem]{Definition}
\newtheorem{remark}[theorem]{Remark}

\newdimen\pIR
\newcommand\StevesR{{\rm I\kern\pIR R}}

\def\bvec#1{{\mbox{\boldmath $#1$}}}

\def\prob#1#2{{\mathbb{P}}_{#1}\left[ #2 \right]}

\newcommand\symExpec{\operatorname{\mathbb{E}}\displaylimits}

\def\expec#1#2{\symExpec_{#1}  #2 }

\def\defeq{\stackrel{\mathrm{def}}{=}}
\def\setof#1{\left\{#1  \right\}}

\def\trace#1{\mathrm{Tr} \left(#1 \right)}

\def\abs#1{\left|#1  \right|}

\def\norm#1{\left\| #1 \right\|}

\def\Complex#1{\mathbb{C}^{#1}}

\newcommand{\imag}[1]{\mathsf{Im} (#1)}

\def\Branden{Br\"{a}nd\'{e}n}

\newcommand{\C}{\mathbb{C}}

\newcommand{\R}{\mathbb{R}}

\newcommand{\mydet}[1]{\det\left(#1\right)}
\newcommand{\charp}[2]{\chi \left[#1 \right] \left(#2 \right)}
\newcommand{\mixed}[2]{\mu \left[#1 \right] \left(#2 \right)}

%% file: ksintro.tex
\section{Introduction}\label{sec:outline}

In their 1959 paper, R. Kadison and I. Singer \cite{orig_ks} posed the following
  fundamental question.
\begin{question}[Kadison-Singer Problem] \label{question:ks} Does every pure state on the (abelian)
von Neumann algebra $\mathbb{D}$ of bounded diagonal operators on
$\ell_2$ have a unique extension to a pure state on $B(\ell_2)$, the
von Neumann algebra of all bounded operators on $\ell_2$?
\end{question}

A positive answer to Question \ref{question:ks} has been shown to be equivalent to a number of conjectures spanning numerous fields, 
  including Anderson's paving conjectures \cite{anderson1979a,anderson1979b,anderson1981},
  Weaver's discrepancy theoretic $KS_{r}$ and $KS_{r}'$ conjectures \cite{weaver},
  the Bourgain-Tzafriri Conjecture \cite{bourgainTzafririUrbana,casazza2006kadison}, and
  the Feichtinger Conjecture and the $R_{\epsilon}$-Conjecture \cite{casazza2005frames};
  and was known to be implied by 
  Akemann and Anderson's projection paving conjecture \cite[Conjecture 7.1.3]{akemann1991lyapunov}.
Many approaches to these problems have been proposed; and, under slightly stronger hypothesis, partial
 solutions have been found by 
  Berman et al. \cite{berman1988matrix}, 
  Bourgain and Tzafriri \cite{BTKS,bourgainTzafririUrbana},
  Paulsen \cite{Paulsen2008},
  Baranov and Dyakonov \cite{BarDya2010},
  Lawton \cite{Lawton2010}, Akemann et al. \cite{akemann2012kadison}, and
  Popa \cite{Popa2013}.
For a discussion of the history and a host of other related conjectures, we refer the reader to \cite{casazza2006detailed}.

We prove these conjectures by proving Weaver's \cite{weaver} conjecture $KS_{2}$ which, as amended by \cite[Theorem 2]{weaver}, says
\begin{conjecture}[$KS_{2}$]\label{conj:weaver}
There exist universal constants $\eta   \geq 2$ and $\theta  > 0$ so that the following holds.
Let $w_{1}, \dots , w_{m} \in \Complex{d}$ satisfy $\norm{w_{i}} \leq 1$ for all $i$
  and suppose
\begin{equation}\label{eqn:weaverCond}
  \sum_{i=1}^m \abs{\langle u, w_{i} \rangle}^{2} = \eta 
\end{equation}
for every unit vector $u \in \Complex{d}$. 
Then there exists a partition $S_{1}, S_{2}$ of $\setof{1, \dots , m}$
  so that
\begin{equation}\label{eqn:weaverNeed}
  \sum_{i \in S_{j}} \abs{\langle u, w_{i} \rangle}^{2} \leq  \eta  - \theta  ,
\end{equation}
for every unit vector $u \in \Complex{d}$ and each $j \in \setof{1,2}$.
\end{conjecture}
Akemann and Anderson's projection paving conjecture \cite[Conjecture 7.1.3]{akemann1991lyapunov}
  follows directly from $KS_{2}$ (see \cite[p. 229]{weaver}).

We also give a proof of Anderson's original paving conjecture, which says
\begin{conjecture}[Anderson Paving]\label{conj:paving} For every $\epsilon>0$ there is an
$r\in\mathbb{N}$ such that for every $n\times n$ self-adjoint complex matrix $T$ with zero
diagonal, there are diagonal projections $P_1,\ldots,P_r$ with $\sum_{i=1}^r P_i=I$ such that
$$ \|P_iTP_i\|\le \epsilon\|T\|\qquad\textrm{for $i=1,\ldots,r$}.$$
\end{conjecture}
A similar conjecture is made by Bourgain and Tzafri \cite[Conjecture 2.8]{bourgainTzafririUrbana}.
One difference between the paving conjecture and 
  $KS_2$ is that the paving conjecture
  can be extended to infinite operators $T\in B(\ell_2)$ by an elementary compactness
  argument \cite{anderson1979a}, which then gives an immediate resolution of
  Kadison-Singer in a manner described in the original paper
  \cite[Lemma 5]{orig_ks}.
On the other hand, the reduction from Kadison-Singer to Akemann and Anderson's
  projection paving conjecture requires non-elementary operator theory.

Our main result follows.  Its proof appears at the end of Section~\ref{sec:barrier}.

\begin{theorem}\label{thm:general}
If $\epsilon > 0$ and
  $v_1, \dots, v_m$ are independent random vectors in $\C^d$ with finite support such that 
\begin{equation}\label{eqn:mainSum}
\sum_{i=1}^m \expec{}{v_{i} v_{i}^{*}} = I_d,
\end{equation}
and
\begin{equation}\label{eqn:mainTrace}
\expec{}{ \norm{v_{i}}^{2}} \leq \epsilon, \text{for all $i$,}
\end{equation}
then 
\[
\prob{}{ \norm{\sum_{i=1}^m v_i v_i^{*}} \leq (1 + \sqrt{\epsilon})^2 } > 0
\]
\end{theorem}

The above theorem may be compared to the concentration inequalities of 
  Rudelson \cite{RudelsonIsotropic} and Ahlswede and Winter \cite{AhlswedeWinter}, which imply in
  our setting that $\|\sum_{i=1}^mv_iv_i^*\|\le C(\epsilon)\cdot\log n$ with high
  probability.
Here we are able to control the deviation at the much smaller
  scale $(1+\sqrt{\epsilon})^2$, but only with nonzero probability.

Our theorem easily implies the following generalization 
  of Conjecture~\ref{conj:weaver}.
\begin{corollary}\label{cor:partition}
Let $r$ be a positive integer and let $u_1, \dots, u_m \in \C^d$ be vectors such that 
$$\sum_{i=1}^m u_iu_i^*=I,$$
and
$\norm{u_i}^2 \leq \delta$ for all $i$.
Then there exists a partition $\{ S_1, \dots S_r \}$ of $[m]$ such that 
\begin{equation}\label{eqn:weaverStrong}
\norm{ \sum_{i \in S_j} {u_i} {u_i}^{*}} \leq
\left(\frac{1}{\sqrt{r}}+\sqrt{\delta}\right)^{2}\qquad\textrm{for
$j=1,\ldots,r$}.
\end{equation}
\end{corollary}

If we set $r=2$ and $\delta = 1/18$, this implies Conjecture~\ref{conj:weaver}
  for  $\eta = 18$ and $\theta = 2$.
To see this,
  set $u_{i} = w_{i} / \sqrt{\eta}$.
Weaver's condition \eqref{eqn:weaverCond} becomes 
  $\sum_{i} u_{i} u_{i}^{*} = I$, and $\delta = 1/\eta$.
When we multiply back by $\eta$,
  the result \eqref{eqn:weaverStrong} becomes \eqref{eqn:weaverNeed}
  with $\eta - \theta = 16$.

Corollary \ref{cor:partition} also implies Conjecture \ref{conj:paving} with
  $r= (6/\epsilon)^{4}$; we defer
  the (slightly more involved) proof to Section \ref{sec:paving}.

\begin{proof}[Proof of Corollary \ref{cor:partition}]
For each $i \in [m]$ and $k \in [r]$, define $w_{i,k} \in \C^{rd}$ to be the direct sum of $r$ vectors from $\C^d$, all of which are $0^d$ (the $0$-vector in $\C^d$) except for the $k^{th}$ one which is a copy of $u_i$.
That is,
\[
  w_{i,1} = \begin{pmatrix}
  u_{i} \\
  0^{d} \\
  \vdots \\
  0^{d}
  \end{pmatrix},
  w_{i,2} = \begin{pmatrix}
  0^{d} \\
  u_{i} \\
  \vdots \\
  0^{d}
  \end{pmatrix},
\text{ and so on}.
\]
Now let $v_1, \dots, v_m$ be independent random vectors such that $v_i$ takes the values $\{ \sqrt{r} w_{i,k} \}_{k=1}^r$ each with probability $1/r$.

These vectors satisfy
\[
  \expec{}{v_{i} v_{i}^{*}}
= 
\begin{pmatrix}
u_{i} u_{i}^{*} & 0_{d \times d} & \dots & 0_{d \times d} \\
0_{d \times d} & u_{i} u_{i}^{*} & \dots & 0_{d \times d} \\
\vdots &  & \ddots &  \vdots \\
0_{d \times d} & 0_{d \times d} & \dots & u_{i} u_{i}^{*},
\end{pmatrix}
\quad \text{and} \quad 
\norm{v_{i}}^{2} = r \norm{u_{i}}^{2} \leq r \delta.
\]
So,
\[
\sum_{i=1}^m  \expec{}{v_{i} v_{i}^{*}} = I_{rd}
\]
and we can apply Theorem~\ref{thm:general} with $\epsilon  = r \delta$ to show that 
  there exists an assignment of each $v_i$ so that  
\[
(1+\sqrt{r \delta})^{2}
\geq 
\norm{\sum_{i=1}^{m} v_{i} v_{i}^{*}}
= 
\norm{
\sum_{k=1}^r
\sum_{i : v_i = w_{i,k}}
\left(\sqrt{r} w_{i,k} \right)
\left(\sqrt{r} w_{i,k} \right)^{*}
}
.
\]
Setting $S_{k} = \{ i : v_i = w_{i,k} \}$, we obtain
\[
\norm{\sum_{i \in S_k} {u_i} {u_i}^{*}} 
= 
\norm{\sum_{i \in S_k} {w_{i,k}} {w_{i,k}}^{*}} 
\leq  
\frac{1}{r}
\norm{
\sum_{k=1}^r
\sum_{i : v_i = w_{i,k}}
\left(\sqrt{r} w_{i,k} \right)
\left(\sqrt{r} w_{i,k}\right)^{*}
}
\leq \left(\frac{1}{\sqrt{r}}+\sqrt{\delta}\right)^{2}.
\]
and this is true for all $k$.
\end{proof}

\section{Overview}

We prove Theorem~\ref{thm:general} using the ``method of interlacing families of polynomials'' introduced in \cite{lifts_paper},
  which we review in Section~\ref{sec:interlacing}.
Interlacing families of polynomials have the property that
  they always contain at least one polynomial whose largest root is at most the
  largest root of the sum of the polynomials in the family.
In Section~\ref{sec:mixed}, we prove that the characteristic polynomials of the matrices
  that arise in Theorem~\ref{thm:general} are such a family.

This proof requires us to consider the expected characteristic polynomials
  of certain sums of independent rank-1 positive semidefinite Hermitian matrices.
We call such an expected polynomial a \textit{mixed characteristic polynomial}.
To prove that the polynomials that arise in our proof are an interlacing family,
  we show that all mixed characteristic polynomials are real rooted.
Inspired by 
  Borcea and \Branden's proof of Johnson's Conjecture~\cite{BBjohnson},
  we do this by constructing multivariate real stable polynomials,
  and then applying operators that preserve real stability
  until we obtain the (univariate) mixed characteristic polynomials.

We then need to bound the largest root of the expected characteristic polynomial.
We do this in Section~\ref{sec:barrier} through a multivariate generalization of the barrier function argument
  of Batson, Spielman, and Srivastava~\cite{BSS}.
The original argument essentially considers the behavior of the roots
  of a real rooted univariate polynomial $p(x)$ under the operator
  $1 - \partial / \partial_{x}$.
It does this by keeping track of an upper bound on the roots of the polynomial,
  along with a measure of how far above the roots this upper bound is.
We refer to this measure as the ``barrier function''.

In our multivariate generalization, we consider a vector $x$ to be {\em above
  the roots} of a real stable multivariate polynomial $p(x_1,\ldots,x_m)$ if
  $p(y_{1}, \dots , y_{m})$ is non-zero
  for every vector $y$ that is at least as big as $x$ in every coordinate.
The value of our multivariate barrier function at $x$ is the vector of the univariate
  barrier functions obtained by restricting to each coordinate.
We then show that we are able to control the values of the barrier function when operators
  of the form 
  $1 - \partial / \partial_{x_{i}}$
  are applied to the polynomial.  
Our proof is inspired by a method used by Gurvits~\cite{gurvitsOne}
  to prove the van der Waerden Conjecture and a generalization 
  by Bapat \cite{bapat1989mixed} of
  this conjecture to mixed discriminants.
Gurvits's proof examines a sequence of polynomials similar to those we construct in our proof,
  and amounts to proving a lower bound on the constant term of the mixed characteristic polynomial.

%% file: kspreliminaries.tex
\section{Preliminaries}\label{sec:prelim}

For an integer $m$, we let $[m] = \setof{1, \dots , m}$.
We write $\binom{[m]}{k}$ to indicate the collection of subsets of $[m]$ having $k$ elements.
When $z_{1}, \dots , z_{m}$ are variables and $S \subseteq [m]$, we define
  $z^{S} = \prod_{i \in S} z_{i}$.

We write $\partial_{z_{i}}$ to indicate the operator that performs partial differentiation
  in $z_{i}$, $\partial / \partial z_{i}$.
For a multivariate polynomial $p (z_{1}, \dots , z_{m})$
  and a number $x$,
  we write $p (z_{1}, \dots , z_{m}) \big|_{z_{1} = x}$
  to indicate the restricted polynomial in $z_{2}, \dots , z_{m}$
  obtained by setting $z_{1}$ to $x$. 

As usual, we write $\norm{x}$ to indicate the Euclidean $2$-norm of a vector $x$.
For a matrix $M$, we indicate the operator norm by $\norm{M} = \max_{\norm{x} = 1} \norm{Mx}$.
When $M$ is Hermitian 
positive semidefinite, 
we recall that this is the largest eigenvalue of $M$.

We write $\mathbb{P}$ and $\mathbb{E}$ for the probability of an event and for the
  expectation of a random variable, respectively.

\subsection{Interlacing Families}\label{sec:interlacing}

We now recall the definition of interlacing families of polynomials from \cite{lifts_paper},
  and its main consequence.
We say that a univariate polynomial is \textit{real rooted} if all of its coefficients and roots are real.

\begin{definition}\label{def:interlacing}
We say that a real rooted
  polynomial $g(x) = \alpha_{0} \prod_{i=1}^{n-1} (x - \alpha_{i})$ \emph{interlaces} a
  real rooted polynomial 
  $f(x) = \beta_{0} \prod_{i=1}^{n} (x - \beta_{i})$ if
\[
  \beta_{1} \leq \alpha_{1} \leq \beta_{2} \leq \alpha_{2} \leq \dotsb \leq
  \alpha_{n-1}\leq \beta_{n}.
\]
We say that polynomials $f_{1} , \dotsc , f_{k} $ have a \emph{common interlacing}
  if there is a polynomial $g $ so that $g $ interlaces $f_i $ for each $i$.
\end{definition}

In \cite{lifts_paper}, we proved the following elementary lemma that shows the utility of having a common interlacing.

\begin{lemma}\label{lem:interlacing}
Let $f_{1}, \dotsc , f_{k}$ be polynomials of the same degree that are real-rooted and have positive leading coefficients.
Define
\[
  f_{\emptyset} = \sum_{i=1}^{k} f_{i}.
\]
If $f_{1}, \dotsc , f_{k}$ have a common interlacing,
  then
  there exists an $i$ so that the largest root of $f_{i}$ is at most the largest root of
  $f_{\emptyset}$.
\end{lemma}

In many cases of interest, we are faced with polynomials that are indexed
  naturally by a cartesian product, and it is beneficial to apply Lemma~\ref{lem:interlacing}
  inductively to subcollections of the polynomials rather than at once.
This inspires the following definition from \cite{lifts_paper}:

\begin{definition}\label{def:family}
Let $S_{1}, \dots , S_{m}$ be finite sets and for every assignment $s_{1}, \dots, s_{m} \in S_{1} \times \dots \times S_{m}$
  let $f_{s_{1}, \dots , s_{m}} (x)$ be a real-rooted degree $n$ polynomial with positive leading coefficient.
For a partial assignment $s_1, \dots, s_k \in S_1 \times \ldots \times S_k$ with $k < m$, define
\[
  f_{s_{1},\dots , s_{k}} \defeq
\sum_{s_{k+1} \in S_{k+1}, \dots , s_{m} \in S_{m}}
  f_{s_{1}, \dots ,s_{k}, s_{k+1}, \dots , s_{m}},
\]
as well as
\[
  f_{\emptyset} \defeq \sum_{s_{1} \in S_{1}, \dots , s_{m} \in S_{m}}
  f_{s_{1}, \dots , s_{m}}.
\]

We say that the polynomials $\{f_{s_{1}, \dots , s_{m}} \}$
  form an \textit{interlacing family} if for all $k=0,\ldots, m-1$, and all
  $s_{1}, \dots, s_{k} \in S_{1} \times \dots \times S_{k}$,
  the polynomials
\[
  \{f_{s_{1}, \dots , s_{k},t}\}_{t\in S_{k+1}}
\]
have a common interlacing.
\end{definition}

\begin{theorem}\label{thm:interlacing}
Let $S_{1}, \dots , S_{m}$
  be finite sets and let
 $\setof{f_{s_{1}, \dots , s_{m}} }$ be an interlacing family of polynomials.
Then, there exists some $s_{1},\dots , s_{m} \in S_{1} \times \dots \times S_{m}$
  so that the largest root of
  $f_{s_{1}, \dots , s_{m}}$ is at most the largest root of $f_{\emptyset}$.
\end{theorem}
\begin{proof}
From the definition of an interlacing family, we know that the polynomials $\{ f_t\}$ for ${t \in S_1}$
  have a common interlacing and that their sum is $f_{\emptyset}$.
So, Lemma~\ref{lem:interlacing} tells us that one of the polynomials $f_{t}$ has largest root at most the largest root of $f_{\emptyset}$.
We now proceed inductively.
For any $s_{1}, \dots , s_{k}$, we know that the polynomials
  $\{f_{s_{1}, \dots , s_k, t}\}$ for $t \in S_{k+1}$ have a common interlacing and
  that their sum is $f_{s_{1}, \dots , s_{k}}$.
So, for some choice of $t$ (say $s_{k+1}$) the largest root of the polynomial
  $f_{s_{1}, \dots , s_{k+1}}$
  is at most the largest root of
  $f_{s_{1}, \dots , s_{k}}$.
\end{proof}

We will prove that the polynomials $\setof{f_{s}}$ defined in Section~\ref{sec:mixed} form an interlacing family.
According to Definition~\ref{def:family}, this requires establishing the
  existence of certain common interlacings.
There is a systematic way to show that polynomials have common interlacings by proving that
  convex combinations of those polynomials are real rooted.
In particular the following result seems to have been discovered a number of times.
It appears as Theorem~$2.1$ of Dedieu~\cite{Dedieu}, (essentially) as
Theorem~$2'$ of Fell~\cite{Fell}, and as (a special case of) Theorem 3.6 of Chudnovsky and
Seymour~\cite{ChudnovskySeymour}.

\begin{lemma}\label{lem:Fisk}
Let $f_1,\ldots, f_k$ be (univariate) polynomials of the same degree with
positive leading coefficients. Then $f_1,\ldots, f_k$ have a common interlacing if
and only if $\sum_{i=1}^k\lambda_i f_i$ is real rooted for all convex
combinations $\lambda_i\ge 0, \sum_{i=1}^k\lambda_i=1$.
\end{lemma}

\subsection{Stable Polynomials}
Our results employ tools from the theory of stable polynomials, a generalization of real rootedness to multivariate polynomials.
For a complex number $z$, let $\imag{z}$ denote its imaginary part.
We recall that a polynomial $p (z_{1}, \dots , z_{m}) \in \Complex{}[z_{1}, \dots , z_{m}]$
  is \textit{stable} if whenever $\imag{z_{i}} > 0$ for all $i$,
  $p (z_{1}, \dots , z_{m}) \not = 0$.
A polynomial $p$ is \textit{real stable} if it is stable and all of its coefficients are real.
A univariate polynomial is real stable if and only if it is real rooted (as defined at the beginning of Section~\ref{sec:interlacing}).

To prove that the polynomials we construct in this paper are real stable,
  we begin with an observation of Borcea and \Branden~\cite[Proposition
  2.4]{BBjohnson}.
\begin{proposition}\label{pro:BBdet}
If $A_{1}, \dots , A_{m}$
  are positive semidefinite Hermitian matrices, then the polynomial
\[
  \mydet{\sum_{i} z_{i} A_{i}}
\]
is real stable.
\end{proposition}

We will generate new real stable polynomials from the one above by applying operators of the form
  $(1 - \partial_{z_{i}})$.
One can use general results, such as Theorem~1.3 of \cite{BBWeylAlgebra}
  or Proposition~2.2 of \cite{LiebSokal}, to prove that these operators
  preserve real stability.
It is also easy to prove it directly using the fact that the analogous operator
  on univariate polynomials preserves stability of polynomials with complex coefficients.
For example, the following theorem appears as Corollary 18.2a in Marden~\cite{Marden},
  and is similar to Corollary 5.4.1 of Rahman and Schmeisser \cite{rahmanSchmeisser}.

\begin{theorem}\label{thm:marden}
If all the zeros of a degree $d$ polynomial $q (z)$ lie in a (closed) circular region $A$,
  then for $\lambda \in \Complex{}$, all the zeros of
\[
  q (z) - \lambda q' (z)
\]
lie in the convex region swept out by translating $A$ in the magnitude and direction
  of the vector $d \lambda$.
\end{theorem}

\begin{corollary}\label{cor:partialRealStable}
If $p \in \R[z_{1}, \dots , z_{m}]$ is real stable,
  then so is
\[
  (1 - \partial_{z_{1}}) p (z_{1}, \dots , z_{m}).
\]
\end{corollary}
\begin{proof}
Let $x_{2}, \dots , x_{m}$ be numbers with positive imaginary part.
Then, the univariate polynomial
\[
  q (z_{1}) =  p (z_{1}, z_{2}, \dots , z_{m}) \big|_{z_{2} = x_{2}, \dots , z_{m} = x_{m}}
\]
is stable.
That is, all of its zeros lie in the circular region consisting of numbers with non-positive
  imaginary part.
As this region is invariant under translation by $d$, $(1-\partial_{z_{1}}) q (z)$
  is stable.
This implies that $(1-\partial_{z_{1}}) p$ has no roots in which all of the variables have positive
  imaginary part.
\end{proof}

We will also use the fact that real stability is preserved under setting
  variables to real numbers (see, for instance, \cite[Lemma 2.4(d)]{wagner}).
\begin{proposition} \label{prop:setreal}
If $p\in\R[z_1,\ldots,z_m]$ is real stable and $a\in \R$, then
$p|_{z_1=a}=p(a,z_2,\ldots,z_m)\in\R[z_2,\ldots,z_m]$ is real stable.
\end{proposition}

\subsection{Facts from Linear Algebra}

For a matrix $M \in \Complex{d \times d}$ we write the characteristic polynomial of $M$ 
  in a variable $x$ as
\[
  \charp{M}{x} = \mydet{xI - M}.
\]



The following identity is sometimes known as the {\em matrix determinant lemma} or 
  the {\em rank-$1$ update formula}.

\begin{lemma}\label{lem:rank1update}
If $A$ is an invertible matrix and $u,v$ are vectors, then
\[
\mydet{A + uv^*} = \mydet{A} (1 + v^* A^{-1} u)
\]
\end{lemma}


We will utilize Jacobi's formula for the derivative of the determinant of a matrix,
  which can be derived from Lemma~\ref{lem:rank1update}.
\begin{theorem}\label{thm:jacobi}
For an invertible matrix $A$ and another matrix $B$ of the same dimensions,
\[
\partial_{t} \mydet{A + tB} =  \mydet{A}\trace{A^{-1} B} .
\]
\end{theorem}

We require two standard facts about traces.
The first is that for a $k$-by-$n$ matrix $A$ and an $n$-by-$k$ matrix $B$,
\[
  \trace{A B} = \trace{B A}.
\]
The second is
\begin{lemma}\label{lem:traceProd}
If $A$ and $B$ are positive semidefinite matrices of the same dimension, then
\[
  \trace{A B} \geq 0.
\]
\end{lemma}
One can prove this by decomposing $A$ and $B$ into sums of rank-1 positive semidefinite matrices,
  using linearity of the trace, and then the first fact about traces.


%


%% file: ksmixed.tex
\section{The Mixed Characteristic Polynomial}\label{sec:mixed}
\begin{theorem}\label{thm:mixed}
Let $v_{1}, \dots , v_{m}$ be independent random column vectors in $\Complex{d}$ with finite support.
For each $i$, let $A_{i} =  \expec{}{v_{i} v_{i}^{*}}$.
Then,
\begin{equation}\label{eqn:mixed1}
\expec{}{\charp{\sum_{i=1}^{m} v_{i} v_{i}^{*}}{x}}
= 
\left(\prod_{i=1}^{m} 1 - \partial_{z_{i}} \right) 
\mydet{x I + \sum_{i=1}^{m} z_{i} A_{i}}
\Big|_{z_{1} = \dots = z_{m} = 0}.
\end{equation}
\end{theorem}
In particular, the expected characteristic polynomial of a sum of independent
  rank one Hermitian matrices is a function of the covariance matrices $A_i$. 
We call this polynomial the \textit{mixed characteristic polynomial} of
  $A_{1}, \dots , A_{m}$, and denote it by $\mixed{A_{1}, \dots , A_{m}}{x}$.

The proof of Theorem~\ref{thm:mixed} relies on the following simple identity
  which shows that random rank one updates of determinants correspond in a natural way to
  differential operators.

\begin{lemma}\label{cor:jameslee}
For every square matrix $A$ and random vector $v$, we have
\begin{equation}\label{eqn:jameslee}
\expec{}{\mydet{A - vv^*}} 
= (1 - \partial_t) \mydet{A + t \expec{}{vv^*}}\big|_{t=0}
\end{equation}
\end{lemma}
\begin{proof}
First,
  assume $A$ is invertible.
By Lemma~\ref{lem:rank1update}, we have
\begin{align*}
\expec{}{\mydet{A - vv^*}}
&= \expec{}{\mydet{A}(1 - v^*A^{-1}v)}
\\&= \expec{}{\mydet{A}(1 - \trace{A^{-1}vv^*})}
\\&= \mydet{A} - \mydet{A}\expec{}{\trace{A^{-1}vv^*}}
\\&= \mydet{A} - \mydet{A}{\trace{A^{-1} \expec{} vv^*}}
\end{align*}
On the other hand, by Theorem~\ref{thm:jacobi}, we have
\[
(1 - \partial_t) \mydet{A + t \expec{}{vv^*}}
= \mydet{A + t \expec{}{vv^*}} - \mydet{A}\trace{A^{-1}\expec{}{vv^*}}.
\]
The claim follows by setting $t = 0$.

If $A$ is not invertible, we can choose a sequence of invertible matrices
  that approach $A$.
Since identity \eqref{eqn:jameslee} holds for each
  matrix in the sequence and the two sides are polynomials in the entries of the matrix, 
  a continuity argument implies that the identity must hold
  for $A$ as well.
\end{proof}

We prove 
 Theorem~\ref{thm:mixed} by applying this lemma inductively.

\begin{proof}[Proof of Theorem~\ref{thm:mixed}]
We will show by induction on $k$ that  for every matrix $M$,
\[
\expec{}{\mydet{M - \sum_{i=1}^{k} v_{i} v_{i}^{*}}}
=\left(\prod_{i=1}^{k} 1 - \partial_{z_{i}} \right) 
\mydet{M + \sum_{i=1}^{k} z_{i} A_{i}}
\Big|_{z_{1} = \dots = z_{k} = 0}.
\]
The base case $k=0$ is trivial. Assuming the claim holds for $i<k$, we have:
\begin{align*}
\expec{}{\mydet{M - \sum_{i=1}^{k} v_{i} v_{i}^{*}}}
&= \expec{v_1, \dots, v_{k-1}}{\expec{v_k}{\mydet{M - \sum_{i=1}^{k-1}
v_{i} v_{i}^{*} - v_{k} v_{k}^{*}}}}\quad\textrm{by independence}
\\&= \expec{v_1, \dots, v_{k-1}}{(1 - \partial_{z_k})\mydet{M - \sum_{i=1}^{k-1} v_{i} v_{i}^{*} + z_k A_k}}\big|_{z_k = 0}
\qquad\textrm{by Lemma \ref{cor:jameslee}}
\\&= (1 - \partial_{z_k})\expec{v_1, \dots, v_{k-1}}{\mydet{M+z_kA_k - \sum_{i=1}^{k-1} v_{i} v_{i}^{*} }}\big|_{z_k = 0}
\qquad\textrm{by linearity}
\\&= (1 - \partial_{z_k})\left(\prod_{i=1}^{k-1} 1 - \partial_{z_{i}} \right)
\mydet{M+z_kA_k + \sum_{i=1}^{k-1} z_{i} A_{i} }\big|_{z_1=\dots = z_{k-1} = 0} \big|_{z_k = 0}
\\&=\left(\prod_{i=1}^{k} 1 - \partial_{z_{i}} \right) \mydet{M + \sum_{i=1}^{k}
z_{i} A_{i} }\big|_{z_1=\dots = z_{k} = 0},
\end{align*}
as desired.
\end{proof}

\begin{remark}\label{rmk:jameslee}
The proof of Theorem~\ref{thm:mixed} given here (using induction and Lemma~\ref{cor:jameslee}) was suggested to us by James Lee.
The slightly longer proof that appeared in our original manuscript was not inductive; rather, it utilized the Cauchy--Binet formula to show
the equality of each coefficient.
\end{remark}

Now it is immediate from Proposition~\ref{pro:BBdet} and Corollary~\ref{cor:partialRealStable} 
  that the mixed characteristic polynomial is real rooted.

\begin{corollary}\label{cor:mixedStable}
The mixed characteristic polynomial of positive semidefinite matrices
  is real rooted.
\end{corollary}
\begin{proof}
Proposition~\ref{pro:BBdet} tells us that 
\[
\mydet{xI+\sum_{i=1}^m z_{i} A_{i}}
\]
  is real stable.
Corollary~\ref{cor:partialRealStable} tells us that
\[
  \left(\prod_{i=1}^{m} 1 - \partial_{z_{i}} \right)
  \mydet{xI+\sum_{i=1}^{m} z_{i} A_{i}}
\]
is real stable as well.
Finally, Proposition \ref{prop:setreal} shows that setting all of the $z_i$ to
  zero preserves real stability.
As the resulting polynomial is univariate, it is real rooted.
\end{proof}

Finally, we use the real rootedness of mixed characteristic polynomials to 
show that every sequence of independent finitely supported random
vectors $v_1,\ldots, v_m$ defines an interlacing family. 
Let $l_{i}$ be the size of the support of the random vector $v_{i}$,
  and let $v_{i}$ take the values $w_{i,1}, \dots , w_{i,l_{i}}$ 
  with probabilities
 $p_{i,1}, \dots , p_{i,l_{i}}$.
For $j_{1} \in [l_{1}], \dots , j_{m} \in [l_{m}]$, define
\[
  q_{j_{1}, \dots , j_{m}} =
  \left(\prod_{i=1}^{m} p_{i,j_{i}} \right)
  \charp{\sum_{i=1}^{m} w_{i,j_{i}} w_{i,j_{i}}^{*} }{x}.
\]
\begin{theorem}\label{thm:mixedInterlacing}
The polynomials $q_{j_{1}, \dots , j_{m}}$ form an interlacing family.
\end{theorem}

\begin{proof}
For $1 \leq k \leq m$ and $j_{1} \in [l_{1}], \dots , j_{k} \in [l_{k}]$, define
\[
q_{j_{1}, \dots , j_{k}} (x)
 = 
\left(\prod_{i=1}^{k} p_{i,j_{i}} \right)
\expec{v_{k+1}, \dots , v_{m}}{\charp{\sum_{i=1}^{k} w_{i,j_{i}} w_{i,j_{i}}^{*} + \sum_{i=k+1}^{m} v_{i} v_{i}^{*} }{x}}.
\]
Also, let 
\[
q_{\emptyset} (x) = 
\expec{v_{1}, \dots , v_{m}}{\charp{\sum_{i=1}^{m} v_{i} v_{i}^{*} }{x}}.
\]
We need to prove that for every partial assignment $j_{1}, \dots , j_{k}$
(possibly empty),
  the polynomials
\[
\left\{q_{j_{1}, \dots , j_{k},t} (x)\right\}_{t=1,\ldots,l_{k+1}}
\]
have a common interlacing.

By Lemma~\ref{lem:Fisk}, it suffices to prove that for every nonnegative
$\lambda_1,\ldots,\lambda_{l_{k+1}}$ summing to one,
  the polynomial
\[
\sum_{t=1}^{l_{k+1}}\lambda_t q_{j_{1}, \dots , j_{k},t} (x)
\]
is real rooted.
To show this, let
  $u_{k+1}$ be a random vector that equals
  $w_{k+1,t}$ with probability $\lambda_t.$
Then, the above polynomial equals
\[
\left(\prod_{i=1}^{k} p_{i,j_{i}} \right)
\expec{u_{k+1}, v_{k+2}, \dots , v_{m}}{\charp{\sum_{i=1}^{k} w_{i,j_{i}}
w_{i,j_{i}}^{*} + u_{k+1}u_{k+1}^{*} + \sum_{i=k+2}^{m} v_{i} v_{i}^{*} }{x}},
\]
which is a multiple of a mixed characteristic polynomial and is thus real rooted by
  Corollary~\ref{cor:mixedStable}.
\end{proof}


%% file: ksbarriers.tex
\def\above#1{\mathsf{Ab}_{#1}}
\section{The Multivariate Barrier Argument}\label{sec:barrier}

In this section we will prove an upper bound on the roots of the mixed
 characteristic polynomial $\mixed{A_1, \dots, A_m}{x}$ as a function of the
 $A_i$, in the case of interest $\sum_{i=1}^mA_i = I$.
Our main theorem is:
\begin{theorem}
\label{thm:mixedbound} Suppose
$A_1,\ldots,A_m$ are Hermitian positive semidefinite matrices satisfying $\sum_{i=1}^m A_i
= I$ and $\trace{A_i}\le \epsilon$ for all $i$. Then the largest root of
$\mixed{A_1,\ldots,A_m}{x}$ is at most $(1+\sqrt{\epsilon})^2$.\end{theorem}

We begin by 
  deriving a slightly different expression for 
  $\mixed{A_1,\ldots,A_m}{x}$
  that allows us to reason  separately about the effect of each $A_i$ on
  its roots.
\begin{lemma}\label{lem:mixedAlt}
Let $A_{1}, \dots , A_{m}$ be Hermitian positive semidefinite matrices.
If $\sum_{i} A_{i} = I$,
 then
\begin{equation}\label{eqn:mixed2}
  \mixed{A_{1}, \dots , A_{m}}{x} = 
  \left(\prod_{i=1}^{m} 1 - \partial_{y_{i}} \right)
  \mydet{\sum_{i=1}^{m} y_{i} A_{i}}
  \Big|_{y_{1} = \dots = y_{m} = x}.
\end{equation}
\end{lemma}
\begin{proof}
For any differentiable function $f$, we have
\[
  \partial_{y_{i}} (f (y_{i})) \big|_{y_{i} = z_{i} + x}
 = 
  \partial_{z_{i}} f (z_{i} + x).
\]
So, the lemma follows by substituting
  $y_{i} = z_{i} + x$
  into expression \eqref{eqn:mixed2},
  and observing that it produces the expression
  on the right hand side of \eqref{eqn:mixed1}.
\end{proof}

Let us write
\begin{equation}\label{eqn:plugq} \mixed{A_1,\ldots,A_m}{x} =
Q(x,x,\ldots,x),\end{equation}
 where $Q(y_1,\ldots,y_m)$ is the multivariate polynomial on the right hand side of
 \eqref{eqn:mixed2}.
The bound on the roots of $\mixed{A_1,\ldots,A_m}{x}$ will follow from a
 ``multivariate upper bound'' on the roots of $Q$, defined as follows.
\begin{definition} Let $p(z_1,\ldots,z_m)$ be a multivariate polynomial. We say
that $z\in\R^m$ is {\em above} the roots of $p$ if 
\[
p(z+t)>0\qquad\textrm{for all}\qquad t=(t_1,\ldots,t_m)\in\R^m, t_i\ge 0,
\]
i.e., if $p$ is positive on the nonnegative orthant with origin at $z$.
\end{definition}

We will denote the set of points which are above the roots of $p$ by
$\above{p}$.  To prove Theorem~\ref{thm:mixedbound}, it is sufficient by
\eqref{eqn:plugq} to show that
  $(1+\sqrt{\epsilon})^2\cdot\bvec{1}\in\above{Q}$, where $\bvec{1}$ is the all-ones vector.
We will achieve this by an inductive ``barrier function'' argument.
In particular, we will construct $Q$ by iteratively 
  applying operations of
  the form $(1-\partial_{y_i})$, and we will track the locations of the roots of
  the polynomials that arise in this process by studying the evolution of the 
  functions defined below.
\begin{definition} Given a real stable polynomial $p$ and a
point $z=(z_1,\ldots,z_m)\in\above{p}$, the {\em barrier function of
$p$ in direction $i$ at $z$} is 
\[
\Phi^i_p(z) = \frac{\partial_{z_i}p(z)}{p(z)} = \partial_{z_i}\log p(z).
\]
\end{definition}
\noindent Equivalently, we may define $\Phi^i_p$ by
\begin{equation}\label{eqn:concdef}
\Phi^i_p(z_1,\ldots,z_m) = \frac{q_{z,i}'(z_i)}{q_{z,i}(z_i)} =
\sum_{j=1}^r\frac{1}{z_i-\lambda_j},\end{equation}
where the univariate restriction
\begin{equation}\label{eqn:restrictdef}
q_{z,i}(t) = p(z_1,\ldots,z_{i-1},t,z_{i+1},\ldots,z_m)\end{equation}
has roots $\lambda_1,\ldots,\lambda_r$, which are real by Proposition
\ref{prop:setreal}.

Although the $\Phi^i_p$ are $m-$variate functions, the properties that we
  require of them may be deduced by considering their bivariate restrictions.
We establish these properties by exploiting the following powerful
  characterization of bivariate
  real stable polynomials.
It is stated in the form we want by
 Borcea and \Branden \ \cite[Corollary 6.7]{BBWeylAlgebra}, and is proved
 using an adaptation of a result of Helton and Vinnikov~\cite{HeltonVinnikov}
 by Lewis, Parrilo and Ramana~\cite{lax}.

\begin{lemma}\label{lem:bivariate}
If $p(z_1,z_2)$ is a bivariate real
stable polynomial of degree exactly $d$, then there exist $d$-by-$d$ positive semidefinite matrices $A,B$ and a
Hermitian matrix $C$ such that 
\[
p(z_1,z_2) = \pm\det(z_1 A+z_2 B+C).
\]
\end{lemma}
\begin{remark}\label{rem:bivariate}
We can also conclude that for every $z_{1}, z_{2} > 0$,
 $z_{1} A + z_{2} B$ must be positive definite.
If this were not the case, then there would be a
 nonzero vector that is in the nullspace of both $A$ and $B$.
This would cause the degree of the polynomial to be lower than $d$.
\end{remark}

The two analytic properties of the barrier functions that we use are
  that, above the roots of a polynomial, they are nonincreasing and convex in every coordinate.
\begin{lemma}\label{lem:monotone}
Suppose $p$ is real stable and
$z\in\above{p}$. Then for all $i,j\le m$ and $\delta\ge 0$,
\begin{align}
\label{eqn:mono}\Phi^i_p(z+\delta e_j) &\le\Phi^i_p(z), \text{ and} \qquad&\textrm{(monotonicity)}\\
\label{eqn:conv}\Phi^i_p(z+\delta e_j) &\le \Phi^i_p(z) + \delta\cdot\partial_{z_j}\Phi^i_p(z+\delta e_j) 
&\qquad\textrm{(convexity).}
\end{align}
\end{lemma}
\begin{proof} 
If $i=j$, then we consider the real-rooted univariate restriction $q_{z,i}(z_i) =
\prod_{k=1}^r(z_i-\lambda_k)$ defined in
\eqref{eqn:restrictdef}. Since $z\in\above{p}$ we know that $z_i>\lambda_k$ for all
$k$. Monotonicity follows immediately by
considering each term in \eqref{eqn:concdef}, and convexity is easily
established by computing
\[
\partial^2_{z_i}\left(\frac{1}{z_i-\lambda_k}\right) =
\frac{2}{(z_i-\lambda_k)^3}>0\qquad\textrm{as $z_i>\lambda_k$}.
\]

In the case $i\neq j$ we fix all variables other than $z_i$ and $z_j$ and consider the
bivariate restriction
$$q_{z,ij}(z_i,z_j) = p(z_1,\ldots,z_m).$$
By Lemma \ref{lem:bivariate} there are Hermitian positive semidefinite $B_i,B_j$ and
a Hermitian matrix $C$ such that
$$q_{z,ij}(z_i,z_j)=\pm\det(z_i B_i+z_j B_j+C).$$
Remark~\ref{rem:bivariate} allows us to conclude that the sign is positive:
 for sufficiently large $t$, $t (B_{i} + B_{j}) + C$ is positive definite
 and for $t \geq  \max (z_{1}, z_{2})$ $q_{z,ij} (t,t) > 0$.

The barrier function in direction $i$ can now be simply expressed as
\begin{align*}
\Phi_p^i(z) &=
\frac{\partial_{z_i}\det(z_iB_i+z_jB_j+C)}{\det(z_iB_i+z_jB_j+C)}
\\&=
\frac{\det(z_iB_i+z_jB_j+C)\trace{(z_iB_i+z_jB_j+C)^{-1}B_i}}{\det(z_iB_i+z_jB_j+C)}
\quad\textrm{by Theorem~\ref{thm:jacobi}}
\\&=\trace{(z_iB_i+z_jB_j+C)^{-1}B_i}
\end{align*}
Let $M=(z_iB_i+z_jB_j+C)$. 
As $z \in \above{p}$ and $B_{i} + B_{j}$ is positive definite, we can conclude that $M$ is positive definite:
 if it were not, there would be a $t$ for which
 $\mydet{M + t (B_{i} + B_{j})} = 0$.
We now write
\begin{align*}
  \Phi_p^i(z+\delta e_j) &= \trace{(M+\delta B_{j})^{-1}B_{i}}\\
   &= \trace{M^{-1} (I+\delta B_{j}M^{-1})^{-1}B_{i}}\\
   &= \trace{(I+\delta B_{j}M^{-1})^{-1}B_{i} M^{-1} }.
\end{align*}
For $\delta$ sufficiently small, we may expand $(I + \delta B_{j}M^{-1})^{-1}$
  in a power series as
\[
  I - \delta B_{j} M^{-1}  + \delta^{2} (B_{j} M^{-1} )^{2} +
  \sum_{\nu \geq 3} (- \delta B_{j} M^{-1})^{\nu }.
\]
Thus,
\[
  \partial_{z_{j}} \Phi_p^i(z)
=
  - \trace{B_{j} M^{-1} B_{i} M^{-1}}.
\]
To see that this is non-positive, 
  and thereby prove \eqref{eqn:mono},
  observe that both $B_{j}$ and $M^{-1}B_{i} M^{-1}$
  are positive semidefinite, and recall from Lemma~\ref{lem:traceProd} that the 
  trace of the product of positive semidefinite matrices is non-negative.
To prove convexity, observe that the second derivative is non-negative because
\[
  \partial_{z_{j}}^{2} \Phi_p^i(z)
=
 \trace{(B_{j} M^{-1} )^{2} B_{i} M^{-1}} 
= 
 \trace{(B_{j} M^{-1} B_{j}) (M^{-1} B_{i} M^{-1})}
\]
is also the trace of the product of positive semidefinite matrices.

Inequality \eqref{eqn:conv} is equivalent to convexity in direction $e_j$ and
 may be obtained by observing that $f(x+\delta)\le f(x)+\delta f'(x+\delta)$
 for any convex differentiable $f$.
\end{proof}

\begin{remark} \label{rmk:renegar}
There are other ways of proving Lemma~\ref{lem:monotone} that go through more elementary
  techniques than those used by  Helton and Vinnikov~\cite{HeltonVinnikov}.
James Renegar has pointed out that it follows from
  Corollary 4.6 of \cite{BGLS}.
Terence Tao~\cite{TaoBlog} has also presented a more elementary proof.
\end{remark}

Recall that we are interested in finding points that lie in $\above{Q}$, where
  $Q$ is generated by applying several operators of the form $1-\partial_{z_i}$ to
  the polynomial $\det(\sum_{i=1}^m z_iA_i)$.
The purpose of the ``barrier functions'' $\Phi_p^i$ is to allow us to
  reason about the relationship between $\above{p}$ and
  $\above{p-\partial_{z_i}p}$;
in particular, the monotonicity property alone immediately implies the following
  statement.
\begin{lemma}\label{lem:above} Suppose that $p$ is real stable, that  $z$
 is above its roots, and  that $\Phi_p^i(z)<1$. Then $z$ is above the roots of
$p-\partial_{z_i}p$.\end{lemma}
\begin{proof}
Let $t$ be a nonnegative vector. As $\Phi$ is nonincreasing in
each coordinate we have $\Phi_p^i(z+t)<1$, whence
$$\partial_{z_i}p(z+t)<p(z+t) \implies (p-\partial_{z_i}p)(x+t)>0,$$
as desired.\end{proof}

Lemma~\ref{lem:above} allows us to prove that a vector is above the roots of
  $p - \partial_{z_{i}} p$.
However, it is not strong enough for an inductive argument because the barrier functions
  can increase with each $1-\partial_{z_i}$ operator that we apply.
To remedy this, we will require the barrier functions to be bounded away
  from $1$, and we will compensate for the effect of each $1-\partial_{z_j}$
  operation by shifting our upper bound away from zero in direction $e_j$.
In particular, by exploiting the convexity properties of the $\Phi_p^i$,
  we arrive at the following strengthening of Lemma~\ref{lem:above}.
\begin{lemma}\label{lem:barrier}
Suppose that $p(z_1,\ldots,z_m)$ is real stable, that $z\in\above{p}$, and that $\delta>0$ satisfies
\begin{equation}\label{eqn:deltaCond}
\Phi^j_p(z) \leq 1 - \frac{1}{\delta}.
\end{equation}
Then for all $i$, 
\[
\Phi^i_{p-\partial_{z_j}p}(z+\delta e_j)\leq\Phi^i_p(z).
\]
\end{lemma}
\begin{proof} 
We will write $\partial_i$ instead of $\partial_{z_i}$ to ease notation. 
We begin by computing
an expression for $\Phi^i_{p-\partial_{j}p}$ in terms of
$\Phi_p^j,\Phi_p^i$, and $\partial_{j}\Phi_p^i$:
\begin{align*}
  \Phi^{i}_{p - \partial_{j}p}
& = 
  \frac{\partial_i (p - \partial_{j} p)}{p  -  \partial_{j} p}
\\
& = \frac{\partial_i \left((1-\Phi^{j}_p) p\right)}{(1-\Phi^{j}_p) p}
\\
&= \frac{(1-\Phi^{j}_p) (\partial_i p)}{(1-\Phi^{j}_p) p} 
  + \frac{(\partial_{i} (1-\Phi^{j}_p)) p}{(1-\Phi^{j}_p) p}
\\
&= \Phi^i_p- \frac{\partial_{i} \Phi^{j}_{p} }
   {1 -  \Phi^{j}_{p}}.
\\
&= \Phi^i_p- \frac{\partial_{j} \Phi^{i}_{p} }
   {1 -  \Phi^{j}_{p}},
\end{align*}
as $\partial_{i} \Phi^{j}_{p} = \partial_{j} \Phi^{i}_{p}$.
We would like to show that $\Phi^i_{p-\partial_{j} p}(z+\delta e_{j} )\leq
  \Phi^i_p(z)$. 
By the above identity this is equivalent to
\[
 -\frac{\partial_{j}  \Phi_p^i(z+\delta e_{j} )}{1-\Phi_p^{j}(z+\delta e_{j} )} \le
\Phi^i_p(z) - \Phi^i_p(z+\delta e_{j} ).
\]

By part \eqref{eqn:conv} of Lemma~\ref{lem:monotone}, 
\[
\delta\cdot (-\partial_{j}  \Phi^i_p(z+\delta e_{j} ))\leq \Phi^i_p(z)-
\Phi^i_p(z+\delta e_{j} ).
\]
Thus it is sufficient to establish that
\begin{equation}\label{le1} -\frac{\partial_{j}  \Phi_p^i(z+\delta e_{j} )}{1-\Phi_p^{j}(z+\delta e_{j} )} \le
\delta\cdot (-\partial_{j}  \Phi^i_p(z+\delta e_{j} )).
\end{equation}

From part \eqref{eqn:mono} of Lemma \ref{lem:monotone}, we know that
$(-\partial_{j} \Phi_p^i(z+\delta e_{j} ))\ge 0$; so, we may divide both sides
  of \eqref{le1} 
 by this term
to obtain
\begin{equation}\label{eqn:le3}
\frac{1}{1-\Phi^i_p(z+\delta e_{j} )} \le \delta.
\end{equation}
Applying Lemma \ref{lem:monotone} once more we observe that $\Phi^{j}_p(z+\delta
e_{j} )\le \Phi^{j}_p(z)$, and conclude that \eqref{eqn:le3} is implied by
\[
\frac{1}{1-\Phi^{j}_p(z)}\le \delta,
\]
which is implied by \eqref{eqn:deltaCond}.
\end{proof}

We now have the necessary tools to prove the main theorem of this section.
\begin{proof}[Proof of Theorem~\ref{thm:mixedbound}]
Let 
\[
  P (y_{1}, \dots , y_{m}) = \mydet{\sum_{i=1}^{m} y_{i} A_{i}}.
\]
Set
\[
  t = \sqrt{\epsilon} + \epsilon .
\]
As all of the matrices $A_{i}$ are positive semidefinite
  and
\[
  \mydet{t \sum_{i} A_{i}} = \mydet{t I} > 0,
\]
the vector $t \bvec{1}$ is above the roots of $P$.

By Theorem~\ref{thm:jacobi},
\[
  \Phi^{i}_{P} (y_{1}, \dots , y_{m})
=
  \frac{\partial_{i} P (y_{1}, \dots , y_{m})}{P (y_{1}, \dots , y_{m})}
=
\trace{\left(\sum_{i=1}^{m} y_{i} A_{i} \right)^{-1} A_{i}}.
\]
So,
\[
  \Phi^{i}_{P} (t \bvec{1})
=
  \trace{A_{i}} / t
\leq 
  \epsilon / t
= 
  \epsilon / (\epsilon + \sqrt{\epsilon}),
\]
which we define to be $\phi$.
Set
\[
  \delta = 1/ (1-\phi) = 1 + \sqrt{\epsilon }.
\]

For $k \in [m]$, define
\[
  P_{k} (y_{1}, \dots , y_{m})
=
  \left(\prod_{i=1}^{k} 1 - \partial_{y_{i}} \right) P (y_{1}, \dots , y_{m}).
\]
Note that $P_{m} = Q$.

Set $x^{0}$ to be the all-$t$ vector, and 
  for $k \in [m]$ define
  $x^{k}$ to be the vector that is $t+\delta$
 in the first $k$ coordinates and $t$ in the rest.
By inductively applying Lemmas~\ref{lem:above} and \ref{lem:barrier},
  we prove that $x^{k}$ is above the roots of $P_{k}$,
  and that for all $i$
\[
   \Phi^{i}_{P_{k}} (x^{k}) \leq \phi .
\]

It follows that the largest root of
\[
  \mixed{A_{1}, \dots , A_{m}}{x}
= 
  P_{m} (x, \dots , x)
\]
is at most 
\[
t + \delta =
  1 + \sqrt{\epsilon} + \sqrt{\epsilon} +  \epsilon  = (1+\sqrt{\epsilon})^{2}.
\]
\end{proof}

\begin{proof}[Proof of Theorem~\ref{thm:general}]
Let $A_{i} = \expec{}{v_{i} v_{i}^{*}}$.
We have 
\[
\trace{A_{i}} = \expec{}{\trace{v_{i} v_{i}^{*}}}
  = \expec{}{v_{i}^{*} v_{i}} = \expec{}{\norm{v_{i}}^{2}} \leq \epsilon ,
\]
for all $i$.

The expected characteristic polynomial of the $\sum_{i} v_{i} v_{i}^{*}$
  is the mixed characteristic polynomial $\mixed{A_{1}, \dots , A_{m}}{x}$.
Theorem~\ref{thm:mixedbound} implies that the largest root of this polynomial 
  is at most $(1+\sqrt{\epsilon})^{2}$.

For $i \in [m]$, let $l_{i}$ be the size of the support of the random vector $v_{i}$,
  and let
  $v_{i}$ take the values $w_{i,1}, \dots , w_{i,l_{i}}$ 
  with probabilities
  $p_{i,1}, \dots , p_{i,l_{i}}$.
Theorem~\ref{thm:mixedInterlacing} tells us that the polynomials
  $q_{j_{1}, \dots , j_{m}}$ are an interlacing family.
So, Theorem~\ref{thm:interlacing} implies that there exist
  $j_{1}, \dots , j_{m}$ so that
  the largest root of the characteristic polynomial of
\[
  \sum_{i=1}^{m} w_{i,j_{i}} w_{i,j_{i}}^{*}
\]
is at most $(1 + \sqrt{\epsilon})^{2}$.
\end{proof}

%% file: kspaving.tex
\section{The Paving Conjecture}\label{sec:paving}
The main result of this section is the following quantitative version of
  Conjecture \ref{conj:paving}.
Following \cite{CEKP07}, we will say that a square matrix $T$ can be {\em $(r,\epsilon)$-paved} if there are
  coordinate projections $P_1,\ldots,P_r$ such that $\sum_{i=1}^r P_i=I$ and $\|P_iTP_i\|\le \epsilon\|T\|$
  for all $i$.
\begin{theorem}\label{thm:paving} For every $\epsilon>0$, every zero-diagonal
  complex self-adjoint matrix $T$ can be $(r,\epsilon)-$paved with
  $r=\left(6/\epsilon\right)^4$.
\end{theorem}
To prove this theorem, we rely on the following result of Casazza et al. which says that
  paving arbitrary self-adjoint matrices can be reduced to paving certain projection
  matrices.
Its short proof is based on elementary linear algebra.
\begin{lemma}[Theorem 3 of \cite{CEKP07}]\label{lem:projpaving} Suppose there is a function
$r:\mathbb{R}_+\rightarrow \mathbb{N}$ so that every $2n \times 2n$ projection matrix $Q$ with diagonal
entries equal to $1/2$ can be $(r(\epsilon),\frac{1+\epsilon}{2})$-paved for all
$\epsilon>0$. Then
every $n\times n$ self-adjoint zero-diagonal matrix $T$ can be
$(r^2(\epsilon),\epsilon)$-paved for all $\epsilon>0$.\end{lemma}
\begin{proof}[Proof of Theorem \ref{thm:paving}] Let $Q$ be an arbitrary
$2n\times 2n$ projection matrix with diagonal entries equal to $1/2$. Then
$Q=(u_i^*u_j)_{i,j\in [2n]}$ is the gram matrix of $2n$ vectors $u_1,\ldots,u_{2n}\in\C^n$ with
$\|u_i\|^2=1/2=\delta$. Applying Corollary \ref{cor:partition} to these vectors
for any $r$
 yields a partition $S_1,\ldots,S_r$ of $[2n]$. Letting $P_k$ be the projection onto the
indices in $S_k$, we have for each $k\in [r]$:
\begin{equation}
\norm{P_kQP_k}=\norm{ \left(u_i^*u_j\right)_{i,j\in S_k}}=\norm{ \sum_{i \in S_k} {u_i} {u_i}^{*}} \leq
\left(\frac{1}{\sqrt{r}}+\frac{1}{\sqrt{2}}\right)^{2}
< \frac{1}{2}+\frac{3}{\sqrt{r}}.
\end{equation}
Thus every $Q$ can be $(r,\frac{1+\epsilon}{2})-$paved for $r=36/\epsilon^2$.
Applying Lemma \ref{lem:projpaving} yields Theorem \ref{thm:paving}.

\end{proof}
It is well-known that Theorem \ref{thm:paving} can be extended to arbitrary
  matrices $T$ with 
  zero diagonal at the cost of a further quadratic loss in parameters:
  simply decompose $T=A+iB$ for self-adjoint zero-diagonal matrices $A,B$, and take a product of pavings of $A$ and $B$.

We have not made any attempt to optimize the dependence of $r$ on $\epsilon$ in
  Theorem~\ref{thm:paving}, and leave this as an open question. 
It is known \cite{CEKP07} that it is not possible to do better than
  $r=1/\epsilon^2$.

%% file: ksconclusion.tex
\section{Conclusion}

When $m = d$, the constant coefficient of the mixed characteristic
  polynomial of $A_{1}, \dots , A_{d}$ 
  is the mixed discriminant of $A_{1}, \dots , A_{d}$.
The mixed discriminant has many definitions, among them
\[
  D (A_{1}, \dots , A_{d}) = 
  \left(\prod_{i=1}^{d} \partial_{z_{i}} \right) \mydet{\sum_{i} z_{i} A_{i}}.
\]
See \cite{gurvitsMixed} or \cite{bapatRaghavan}.

When $k < d$, we define
\[
  D (A_{1}, \dots , A_{k})
=
  D (A_{1}, \dots , A_{k}, I, \dots , I) / (d-k)!,
\]
where the identity matrix $I$ is repeated $d-k$ times.
For example $D (A_{1})$ is the just the trace of $A_{1}$.

With this notation, we can write
\[
  \mixed{A_{1}, \dots , A_{m}}{x}
= 
  \sum_{k=0}^{d} x^{d-k} (-1)^{k}
  \sum_{S \in \binom{[m]}{k}} D ((A_{i})_{i \in S}).
\]

When the matrices $A_{1}, \dots , A_{d}$ are diagonal,
  $\mixed{A_{1}, \dots , A_{d}}{x}$ is the matching 
  polynomial defined by Heilmann and Lieb~\cite{heilmannLieb} of the bipartite graph with
  $d$ vertices on each side in which the edge $(i,j)$
  has weight $A_{i} (j,j)$.
When all the matrices have the same trace and their sum is the identity,
  the graph is regular and our bound on the largest root of the mixed
  characteristic polynomial agrees to the first order with that obtained for the matching
  polynomial by Heilmann and Lieb~\cite{heilmannLieb}.
  
We conjecture that among the families of matrices $A_{1}, \dots , A_{m}$
  with $\sum_{i} A_{i} = I$ and $\trace{A_{i}} \leq \epsilon$,
  the largest root of the mixed characteristic polynomial
  is maximized when as many of the matrices as possible equal $\epsilon I / d$,
  another is a smaller multiple of the identity, and the rest are zero.
When all of the matrices have the same trace, $d/m$,
  this produces a scaled associated Laguerre polynomial $L_{d}^{m-d} (m x)$.
The bound that we prove on the largest root of the mixed characteristic
  polynomial agrees asymptotically with the largest root of $L_{d}^{m-d} (mx)$ as $d/m$ is held
  constant and $d$ grows.
Evidence for our conjecture may be found in the work of 
  Gurvits~\cite{gurvitsMixed,gurvitsOne}, who proves that when $m = d$, the constant term of the mixed
  polynomial is minimized when each $A_{i}$ equals $I/d$.

Two natural questions arise from our work.
The first is whether one can design an efficient algorithm
  to find the partitions and pavings that are guaranteed to exist by 
  Corollary~\ref{cor:partition}.
The second is broader.
There are many operations that are known to preserve
  real stability and real rootedness of polynomials (see
  \cite{LiebSokal,BBWeylAlgebra,BBpolyaSchurI,BBpolyaSchurII,pemantle,wagner}).
For a technique like the ``method of characteristic polynomials'' it would be useful to understand what these operations do to the
  roots and the coefficients of the polynomials.


\section{Acknowledgements}
We thank Gil Kalai for making us aware of the Kadison-Singer conjecture
  and for pointing out the resemblance between the paving conjecture and
  the sparsification results of~\cite{BSS}.
We thank Joshua Batson for discussions of this connection at the
  beginning of this research effort.
We thank
  James Lee for suggesting the simple proof of Theorem \ref{thm:mixed}.
Following the initial release of this paper, we received a number of comments and improvements, and we would like to thank everyone that took the time to give us input.